\documentclass[12pt]{amsart}
\usepackage{amsfonts,latexsym}

\newtheorem{theorem}{Theorem}

\newtheorem{thm}[theorem]{Theorem}
\newtheorem{prop}[theorem]{Proposition}

\theoremstyle{definition}

\theoremstyle{remark}
\newtheorem{remark}[theorem]{Remark}

\title[Combinatorial principles, compactness of spaces]
{Combinatorial and model-theoretical principles related to 
regularity of ultrafilters and compactness of topological spaces. I.}

\author[]{Paolo Lipparini} 
\address{Dipartimento di Matematica\\
Viale della Ricerca Scientifica\\
II Universit\`a de Roma (Tor Vergata)\\
I-00133 ROME 
ITALY
}
\urladdr{http://www.mat.uniroma2.it/\textasciitilde lipparin}
\thanks{The author has received support from MPI and GNSAGA.
We wish to expressed our gratitude to X. Caicedo for stimulating discussions and correspondence} 


\begin{document} 

\begin{abstract} 
We begin the study of the consequences of the existence of certain
 infinite matrices. Our present application is to
 compactness of products of topological spaces.
\end{abstract}

\maketitle




Our notation is fairly standard. See, e.~g., \cite{CN,KV,EGT}
for unexplained notation.

Ordinals are denoted by $\alpha, \beta, \gamma,  \dots$
Infinite cardinals are denoted by $\lambda, \mu, \nu, \kappa, \dots$
Inclusion is denoted by $\subseteq$, and  
$\subset$ denotes
strict inclusion. The minus operation between sets is denoted by $\setminus$,
that is,
$X \setminus Y = \{x\in X| x \not \in Y\} $.

We assume the Axiom of Choice.

If
$ (X_ \alpha) _{\alpha < \lambda }  $
are topological spaces, then
$\prod_{\alpha <\lambda } X_ \alpha  $
denotes their product with the Tychonoff topology,
 the smallest
topology under which the canonical
projections are continue maps.

The {\em $ \lambda  $-th power} of a topological space $X$
is the product $\prod_{\alpha <\lambda } X_ \alpha  $,
where $X_ \alpha = X$ for all $ \alpha \in \lambda $.

If $ \kappa  , \lambda  $ are infinite cardinals, a topological space is
said to be  {\em $[ \kappa  , \lambda  ]$-compact}
if and only if
every open cover by at most $ \lambda $  sets  has a subcover by less than 
$\kappa $ sets.

No separation axiom is needed to prove the results of the present paper.

The following characterizations are old and well-known. See \cite[Section 3]{topappl} 
for details, further references and further information about 
$[ \kappa  , \lambda  ]$-compactness.

\begin{prop}\label{regcomp}
For every infinite regular cardinal $ \kappa $ and every topological
space $X$, the following are equivalent. 

(i) $X$ is $[ \kappa  , \kappa   ]$-compact.

(ii) Whenever $(U_ \alpha )_{\alpha < \kappa } $
is a sequence of open sets of $X$, such that 
$ U_ \alpha \subseteq U_{\alpha'}$ for every $ \alpha<\alpha' $,
and such that $\bigcup_{ \alpha < \kappa } U_ \alpha =X$,
then there is an $ \alpha < \kappa $ such that $ U_ \alpha =X$.

(iii) Whenever $(C_ \alpha )_{\alpha < \kappa } $
is a sequence of closed sets of $X$, such that 
$ C_ \alpha \supseteq C_{\alpha'}$ for every $ \alpha<\alpha' $,
and such that $\bigcap_{ \alpha < \kappa } C_ \alpha =\emptyset$,
then there is an $ \alpha < \kappa $ such that $ C_ \alpha =\emptyset$.
 
(iv) For every sequence $(x_ \alpha )_{\alpha < \kappa } $
 of elements of $X$, there exists $x \in X$ such that 
$ |\{ \alpha < \kappa | x_ \alpha \in U \}|=  \kappa $
for every neighbourhood $U$ of $x$. 

(v) (CAP$_\kappa $) Every subset $Y\subseteq X$ with $|Y|=\kappa $ has a complete
accumulation point.
\end{prop}

\begin{thm}\label{lmkprod} 
Suppose that  $ \lambda $, $\mu$  are infinite regular cardinals, 
and $\kappa $ is an infinite cardinal.  
Then the following conditions are equivalent.

\smallskip

(a) There is a family $ (B_{ \alpha , \beta }) _{ \alpha<\mu , \beta<\kappa}  $ 
of subsets of $ \lambda $ such that:

(i) For every $ \beta<\kappa$, $\bigcup _{ \alpha<\mu } B_{ \alpha , \beta  } = \lambda$;

(ii) For every $ \beta<\kappa$ and $ \alpha \leq \alpha ' < \mu  $, 
$ B_{ \alpha , \beta } \subseteq B_{ \alpha' , \beta }$;

(iii) For every function $f : \kappa  \to \mu $ there exists a finite subset
$F \subseteq \kappa  $ such that 
$|\bigcap _{\beta \in F} B_{ f( \beta) , \beta }| < \lambda $.   

\smallskip

(b)
Whenever
$(X_ \beta ) _{ \beta < \kappa }$ is a family of topological spaces
such that  no $X_ \beta $ is 
$[ \mu, \mu]$-compact,
then $X=\prod_{ \beta < \kappa } X_ \beta $
 is not  $[ \lambda , \lambda ]$-compact.

\smallskip

(c) The topological space $ \mu^ \kappa $ is not
$[ \lambda , \lambda ]$-compact, where $ \mu$
is endowed with the topology whose open sets are the
intervals $ [0, \alpha) $ ($ \alpha \leq \mu$), and 
$ \mu^ \kappa $ is endowed with the Tychonoff topology. 
\end{thm} 
 
\begin{remark}\label{future}
In a sequel to this note we shall provide many more 
conditions equivalent to the 
conditions in Theorem \ref{lmkprod}.
The same applies to the conditions we shall introduce
in Theorems \ref{lmkprodbox} and \ref{lmkprodgeneral}.
\end{remark}

\begin{proof}
(a) $\Rightarrow$ (b).
Let $X, (X_ \beta) _{ \beta < \kappa }  $ and
$ (B_{ \alpha , \beta }) _{ \alpha<\mu , \beta<\kappa}  $ 
 be as in the statement
of the theorem.

Since 
no $X_\beta$ is 
$[ \mu, \mu]$-compact, 
and since $ \mu $ is regular, 
by Condition (iv) in Proposition \ref{regcomp}, 
for every $\beta<\kappa$ there is a 
 sequence $\{x_{\alpha , \beta}| \alpha <\mu \} $ 
of elements of $X_ \beta $ 
such that 
every $x \in X_\beta$ has a neighbourhood $U_ \beta$  
in $X_\beta$  such that 
$ |\{ \alpha<\mu | x_{\alpha , \beta} \in U_ \beta  \}|<  \mu  $.

We shall define a sequence 
$(y_\gamma) _{\gamma<\lambda}$ of elements of $ X$ such that 
for every $z \in X$ there is a neighbourhood $U$ in $X$ of $z$
 such that 
$ |\{\gamma<\lambda| y_\gamma \in U\}|<  \lambda$, thus 
 $X$ is not  $[ \lambda, \lambda]$-compact,
again by Condition (iv) in Proposition \ref{regcomp},
and since  $ \lambda $ is supposed to be  regular.

For $\gamma<\lambda$, 
let $y_\gamma =((y_\gamma )_\beta)_{\beta<\kappa} \in \prod_{\beta<\kappa} X_\beta $
be defined by: $(y_\gamma )_\beta =x_{\alpha, \beta}$, where $\alpha $ is the first ordinal
such that $\gamma \in B_{ \alpha , \beta } $ (such an ordinal exists by Condition (i)).

Suppose by contradiction that there is $z \in X$ 
 such that for every neighbourhood $U$ in $X$ of $z$
 $ |\{\gamma<\lambda| y_\gamma \in U\}|=  \lambda$.

Consider the components $(z_\beta)_{\beta<\kappa}$
of $z\in X=\prod_{\beta<\kappa} X_\beta$.
Because of the way we have chosen the 
$ x_{\alpha , \beta}$s,
for each $\beta<\kappa$, $z_\beta$ has a neighbourhood
$U_\beta$ in $X_\beta$ such that 
$ |\{\alpha < \mu| x_{\alpha , \beta}\in U_\beta\}| < \mu  $.
 For every $\beta<\kappa$, fix some $U_\beta$ as above. 
For each $\beta<\kappa$, choose $f(\beta)$ in such a way that
$\mu >f(\beta)>\sup \{\alpha < \mu | x_{\alpha , \beta}\in U_\beta\}$
 (this is possible since $\mu $ is regular, and
$ |\{\alpha < \mu | x_{\alpha , \beta}\in U_\beta\}| < \mu  $).

By Condition (iii) 
there is a finite $F \subseteq \kappa $ such that  
$|\bigcap _{\beta \in F} B_{ f( \beta) , \beta }| < \lambda  $.   
Let $V=\prod_{\beta<\kappa} V_\beta$, where 
$V_\beta=X_\beta$ if $\beta \not \in F$, and
$V_\beta=U_\beta$ if $\beta \in F$.
$V$ is a neighbourhood of $z$ in $X$, since $F$
is finite.

For every $ \gamma<\lambda$ and $\beta<\kappa$,
by  definition,  $(y_\gamma )_\beta = x_{\alpha, \beta}$, for some $\alpha $ 
such that $\gamma \in B_{ \alpha , \beta } $. By the definition of $f$,
if $(y_\gamma )_\beta = x_{\alpha, \beta}\in U_\beta$ then 
$f(\beta )>\alpha $, thus $ \gamma  \in B_{ \alpha , \beta } \subseteq B_{ f( \beta ), \beta  }$,
by Condition (ii).
We have proved that,
for every $\beta<\kappa$, 
$\{\gamma<\lambda| (y_\gamma )_\beta \in U_\beta\} \subseteq B_{ f(\beta), \beta }$.

Thus, by the definition of $V$, we have
$\{\gamma<\lambda| y_\gamma \in V\}=
\bigcap _{\beta \in F} 
\{\gamma<\lambda| (y_\gamma )_\beta \in U_\beta\}
\subseteq \bigcap _{\beta \in F} 
B_{ f(\beta), \beta }$.
 Hence 
$|\{\gamma<\lambda| y_\gamma \in V\}|
\leq |\bigcap _{ \beta  \in F} 
B_{ f(\beta), \beta }| < \lambda $.
This is a contradiction, since 
we have supposed that 
$|\{\gamma<\lambda| y_\gamma \in V\}| =\lambda$,
for every 
neighbourhood $V$ of $z$.  

(b) $\Rightarrow$ (c) is trivial, since $\mu$ is not  
$[ \mu, \mu]$-compact.

(c) $\Rightarrow$ (a).
By Condition (iv) in Proposition \ref{regcomp}
there exists a sequence 
$ (y_ \gamma ) _{ \gamma < \lambda } $
of elements in $ \mu ^ \kappa $
such that   
for every $z \in \mu ^ \kappa $
 there is a neighbourhood $U$ in $\mu ^ \kappa $ of $z$
 such that 
$ |\{\gamma<\lambda| y_\gamma \in U\}|<  \lambda$.

For each $ \gamma < \lambda $,  
$y_\gamma  \in \mu ^ \kappa  $ has the form
$y_\gamma =((y_\gamma )_\beta)_{\beta<\kappa} $.
For $ \alpha < \mu$ and $ \beta < \kappa $ define
$B _{ \alpha, \beta } = 
\{ \gamma < \lambda | (y_\gamma )_\beta  \leq \alpha \}$.
 
Conditions (i) and (ii) in (a) trivially hold.

As for Condition (iii), suppose that $f: \kappa \to \mu$.
Let $z\in \mu^ \kappa $ be defined by 
$z=(f( \beta )) _{ \beta < \kappa } $.
By the first paragraph, 
there is a neighbourhood $U$ in $\mu ^ \kappa $ of $z$
 such that 
$ |\{\gamma<\lambda| y_\gamma \in U\}|<  \lambda$.

Arguing componentwise, this means that
there are a finite set $ F \subseteq \kappa $
and, for each $ \beta \in F$, 
neighbourhoods $U _ \beta $ of $f( \beta )$ in $ \mu$
such that 
$|\bigcap _{\beta \in F} 
\{\gamma<\lambda| (y_\gamma )_\beta \in U_\beta\}|< \lambda $.
Since any
neighbourhood $U _ \beta $ of $f( \beta )$ in $ \mu$
contains $[0,f( \beta )+1)$, we have that 
$(y_\gamma )_\beta \leq f(\beta)$ implies
that $(y_\gamma )_\beta \in U_\beta$. Hence also
$|\bigcap _{\beta \in F} 
\{\gamma<\lambda| (y_\gamma )_\beta \leq f( \beta )\}|< \lambda $.

Thus,
 $\bigcap _{ \beta  \in F} B_{ f(\beta), \beta } =
\bigcap _{ \beta  \in F} 
\{ \gamma < \lambda | (y_\gamma )_\beta  \leq f( \beta ) \}$
 has cardinality $< \lambda $. 
   \end{proof} 

\begin{remark}\label{gen}
In the particular
case $ \lambda = \kappa = \mu^+$  \cite[Lemma 14]{topappl} states 
that Condition (a) in Theorem \ref{lmkprod} is true, and, actually,
we can get $|F|=2$ (the proof elaborates on a variation on a classical
combinatorial device known as an ``Ulam matrix'' \cite{EU}). 
Proposition 15 in \cite{topappl}  then goes on showing that,
in the above particular case $ \lambda = \kappa = \mu^+$,
Condition (b) in Theorem \ref{lmkprod} holds.
 Thus, modulo \cite[Lemma 14]{topappl}, Theorem \ref{lmkprod} generalizes
\cite[Proposition 15]{topappl}. Indeed, our proof of (a) $ \Rightarrow $ (b)
in Theorem \ref{lmkprod} is modelled after the proof of Proposition 15 in \cite{topappl}.
 
 The main results proved in \cite{topappl} had been announced in \cite{abst},
where further results similar to the ones presented here are stated.
 \cite{Cprepr,C,bumi,arxiv,topproc}
also contain related results. We plan to give a unified treatment of all 
these results in a sequel to the present note.
\end{remark}

Theorem \ref{lmkprod} can be generalized for box products.

If  $ \nu$ is a cardinal, and  
$(X_ \beta ) _{ \beta < \kappa }$ is a family of topological spaces, 
then their product can be assigned the 
$\Box^ {<\nu} $ topology, the topology a base of which
is given by all products $(Y_ \beta ) _{ \beta < \kappa }$,
where each $Y_ \beta $ is an open subset of $X_ \beta $,
and $| \{  \beta < \kappa | Y_ \beta \not= X_ \beta \}| < \nu$.   
The product of $(X_ \beta ) _{ \beta < \kappa }$
with the $\Box^ {<\nu} $ topology shall be denoted by 
$\Box^{<\nu}_{ \beta < \kappa } X_ \beta $.

\begin{thm}\label{lmkprodbox} 
Suppose that  $ \lambda $, $\mu$  are infinite regular cardinals, 
and $\kappa $, $ \nu$  are infinite cardinals.  

Then the following conditions are equivalent.

\smallskip

(a) There is a family $ (B_{ \alpha , \beta }) _{ \alpha<\mu , \beta<\kappa}  $ 
of subsets of $ \lambda $ such that:

(i) For every $ \beta<\kappa$, $\bigcup _{ \alpha<\mu } B_{ \alpha , \beta  } = \lambda$;

(ii) For every $ \beta<\kappa$ and $ \alpha \leq \alpha ' < \mu  $, 
$ B_{ \alpha , \beta } \subseteq B_{ \alpha' , \beta }$;

(iii) For every function $f : \kappa  \to \mu $ there exists a subset
$F \subseteq \kappa  $ such that $|F|<\nu$ and 
$|\bigcap _{\beta \in F} B_{ f( \beta) , \beta }| < \lambda $.   

\smallskip

(b)
Whenever
$(X_ \beta ) _{ \beta < \kappa }$ is a family of topological spaces 
such that  no $X_ \beta $ is 
$[ \mu, \mu]$-compact,
then $X=\Box^{<\nu}_{ \beta < \kappa } X_ \beta $
 is not  $[ \lambda , \lambda ]$-compact.

\smallskip

(c) The topological space $ \mu^ \kappa $ is not
$[ \lambda , \lambda ]$-compact, where $ \mu$
is endowed with the topology whose open sets are the
intervals $ [0, \alpha) $ ($ \alpha \leq \mu$), and 
$ \mu^ \kappa $ is endowed with the 
$\Box^{<\nu}$
 topology. 
\end{thm} 

\begin{proof} 
The proof is identical to the proof of Theorem \ref{lmkprod}.
\end{proof} 

Notice that Theorem \ref{lmkprodbox} generalizes
Theorem \ref{lmkprod}, since
the Tychonoff product is just the box product
$\Box^{ <\omega } $. Hence 
Theorem \ref{lmkprod} is the particular case
$\nu=\omega$ of  Theorem \ref{lmkprodbox}.

We have an even more general version of the above theorems.

\begin{thm}\label{lmkprodgeneral} 
Suppose that  $ \lambda $ is an infinite regular cardinal, 
 $\kappa $, $ \nu$  are infinite cardinals,
and $(\mu_ \beta ) _{ \beta < \kappa } $ are infinite regular cardinals.  

Then the following conditions are equivalent.

\smallskip

(a) There is a family $ (B_{ \alpha , \beta }) _{  \beta<\kappa, \alpha<\mu_ \beta}  $ 
of subsets of $ \lambda $ such that:

(i) For every $ \beta<\kappa$, $\bigcup _{ \alpha<\mu_ \beta } B_{ \alpha , \beta  } = \lambda$;

(ii) For every $ \beta<\kappa$ and $ \alpha \leq \alpha ' < \mu_ \beta $, 
$ B_{ \alpha , \beta } \subseteq B_{ \alpha' , \beta }$;

(iii) For every $ f \in \prod _{ \beta < \kappa } \mu _ \beta $ 
there exists a subset
$F \subseteq \kappa  $ such that $|F|<\nu$ and 
$|\bigcap _{\beta \in F} B_{ f( \beta) , \beta }| < \lambda $.   

\smallskip

(b) Whenever
$(X_ \beta ) _{ \beta < \kappa }$ is a family of topological spaces 
such that  for no $ \beta < \kappa $ $X_ \beta $ is 
$[ \mu_ \beta, \mu_ \beta]$-compact,
then $X=\Box^{<\nu}_{ \beta < \kappa } X_ \beta $
is not  $[ \lambda , \lambda ]$-compact.

\smallskip

(c) The topological space $ \Box^{<\nu}_{ \beta < \kappa } \mu_ \beta  $ is not
$[ \lambda , \lambda ]$-compact, where, for each $ \beta < \kappa $,
$ \mu _ \beta $
is endowed with the topology whose open sets are the
intervals $ [0, \alpha) $ ($ \alpha \leq \mu _ \beta $).
\end{thm}

\begin{proof} 
The proof is similar to the proof of Theorem \ref{lmkprod}.
\end{proof}

\end{document}